\documentclass[11pt]{article}

\def\emline#1#2#3#4#5#6{%
       \put(#1,#2){\special{em:moveto}}%
       \put(#4,#5){\special{em:lineto}}}

\usepackage{latexsym}
\usepackage{amsmath}
\usepackage{amssymb}
\usepackage{amsfonts}
\usepackage{amsthm}
\usepackage{stmaryrd}  

\numberwithin{equation}{section}

\makeatletter
\newcommand{\textbiginterleave}{{\fontencoding{U}\fontfamily{stmry}
  \fontseries{m}\fontshape{n}\fontsize{14.4}{\f@baselineskip}
  \selectfont\char57}}
\newcommand{\textbbiginterleave}{{\fontencoding{U}\fontfamily{stmry}
  \fontseries{m}\fontshape{n}\fontsize{17.28}{\f@baselineskip}
  \selectfont\char57}}
\makeatother

\theoremstyle{plain}
\newtheorem{theo}{Theorem}
\newtheorem{lemma}{Lemma}
\newtheorem{cor}{Corollary}
\newtheorem{rem}{Remark}

\newcommand{\noi}{\noindent}
\newcommand{\bn}{\bigskip\noindent}

\newcommand{\bsk}{\bigskip}
\newcommand{\q}{\quad}
\newcommand{\qq}{\qquad}
\newcommand{\e}{\ell}

\newcommand{\lra}{\longrightarrow}
\newcommand{\hra}{\hookrightarrow}
\newcommand{\car}{\curvearrowright}
\newcommand{\rsa}{\rightsquigarrow}
\newcommand{\cale}{\curvearrowleft}
\newcommand{\carle}{\curvearrowright\hspace*{-4.2mm}\curvearrowleft}
\newcommand{\ot}{\otimes}
\newcommand{\C}{{\mathbb C}}
\newcommand{\K}{{\mathbb K}}
\newcommand{\R}{{\mathbb R}}
\newcommand{\N}{{\mathbb N}}
\newcommand{\tsum}{{\textstyle{\sum}}}
\newcommand{\ci}{{\,\scriptstyle{\circ}\,}}

\DeclareMathOperator{\rinf}{inf}

\DeclareMathOperator{\ran}{\rm range}

\newcommand{\bgdv}{\bigg\|}   
\newcommand{\bgl}{\bigg(}
\newcommand{\bgr}{\bigg)}
\newcommand{\bgv}{\bigg|}

\newcommand{\bv}{\big|}
\newcommand{\bdv}{\big\|}    
\newcommand{\bl}{\bigl(}
\newcommand{\br}{\bigr)}
\newcommand{\brb}{\bigr]}
\newcommand{\blb}{\bigl[}

\newcommand{\Bdv}{\Big\|}   
\newcommand{\Bl}{\Bigl(}
\newcommand{\Br}{\Bigr)}

\newcommand{\bs}{\big/}

\newcommand{\ov}{\overline}

\newcommand{\ck}{\mathcal K}

\newcommand{\la}{\langle}
\newcommand{\ra}{\rangle}
\newcommand{\tv}{\interleave}
\newcommand{\btv}{\big|\hspace*{-0.5mm}\big|\hspace*{-0.5mm}\big|}

\newcommand{\blk}{\bigl\{}
\newcommand{\brk}{\bigr\}}
\newcommand{\Blk}{\Bigl\{}
\newcommand{\Brk}{\Bigr\}}
\newcommand{\bglk}{\biggl\{}
\newcommand{\bgrk}{\biggr\}}

\newcommand{\bmv}{\bigm|}

\newcommand{\Bmv}{\Bigm|}

\title{\bf Variants of the Maurey-Rosenthal theorem\\[1ex] for\\[1ex]
quasi K\"othe function spaces}
\author{{\small Andreas Defant}\\
{\small Fachbereich Mathematik}\\
{\small Universit\"at Oldenburg}\\
{\small D-26111 Oldenburg, Germany}\\
{\small E-mail: defant@mathematik.uni-oldenburg.de}}
\date{}

\begin{document}
\maketitle
\begin{abstract}
The Maurey-Rosenthal theorem states that each bounded and linear operator
$T$ from a quasi normed space $E$ into some $L_p(\nu)$ $(0 < p < r <
\infty)$ which satisfies a vector-valued norm inequality
\[
   \bdv\bl\tsum |Tx_k|^r\br^{1/r}\bdv_{L_p} \le
          \bl\tsum\|x_k\|^r_E\br^{1/r}
          \qq\text{for all $x_1,\ldots,x_n\in E$,}
\]
even allows a weighted norm inequality: there is a function $0\le w\in
L_0(\nu)$ such that
\[
   \bgl\int \frac{|Tx|^r}{w}\, d\nu\bgr^{1/r} \le \|x\|_E
     \qq\text{for all $x\in E$.}
\]
Continuing the work of Garcia-Cuerva and Rubio de Francia we give several
scalar and vector-valued variants of this fundamental result within the
framework of quasi K\"othe function spaces $X(\nu)$ over measure spaces.
\end{abstract}


The fundamental Maurey-Rosenthal theorem -- due to Maurey \cite{Mau74}
with
roots in the classical work \cite{Ro73} of Rosenthal -- states that each
bounded and linear operator $T$ from a quasi normed space $E$ into some
$L_p(\nu)$ ($\nu$ a measure and
$0 < p <\infty$) satisfies for $0< p < r < \infty$ a {\em vector-valued
norm inequality}
\begin{equation}\label{eq-0.1}
   \bgdv\bgl \sum^n_{k=1} |Tx_k|^r \bgr^{1/r} \bgdv_{L_p} \le
   \bgl\sum^n_{k=1} \|x_k\|^r_E\bgr^{1/r}
\end{equation}
if and only if it satisfies a {\em weighted norm inequality}
\begin{equation}\label{eq-0.2}
   \bgl\int \frac{|Tx|^r}{w}\, d\nu\bgr^{1/r} \le \|x\|_E\, ;
\end{equation}
more precisely, $T$ fulfills \eqref{eq-0.1} for all choices of vectors
$x_1,\ldots,x_n\in E$ if and only if there exists a weight $0 \le w \in
L_0(\nu)$ (with an appropriate norm estimate) such that \eqref{eq-0.2}
holds for all $x\in E$.

This result is of special interest for operators $T: E\lra F$, $F$ a
subspace of some $L_1(\nu)$. Dually, the famous Grothendieck-Pietsch
domination theorem says that a bounded and linear operator  $T: E\lra F$,
$E$ now a subspace of some $L_\infty(\mu)$ and $F$ a Banach space, allows
for $1\le r < \infty$ a vector-valued norm inequality
\begin{equation}\label{eq-0.3}
   \bgl\sum^n_{k=1} \|Tx_k\|^r_F\bgr^{1/r} \le \bgdv \Bl\sum^n_{k=1}
     |Tx_k|^r\bgr^{1/r}\bgdv_{L_\infty}
\end{equation}
($T$ is $r$-summing) if and only if the weighted norm inequality
\begin{equation}\label{eq-0.4}
   \|Tx\|_F \le \varphi(|x|^r)^{1/r}
\end{equation}
holds true; more accurately: $T$ satisfies \eqref{eq-0.3} for all
$x_1,\ldots, x_n\in E$ iff there is a positive continuous functional
$\varphi\in L_\infty(\mu)'$ (with $\|\varphi\| = 1$) such that
\eqref{eq-0.4} holds for all $x\in E$.

Nowadays such equivalences of vector-valued and weighted norm inequalities
are of fundamental importance in different parts of analysis (see \cite{Ro73},
\cite{Mau74}, \cite{Pie80}, \cite{DJT95}, \cite{Wo91} and \cite{DF93} for
operator theory and geometry in Banach spaces, and \cite{Gar90},
\cite{GR85},
\cite{St93} and \cite{Wo91} for harmonic analysis). In particular, Rubio de
Francia \cite{Ru87}, \cite{Ru85} and later Garcia-Cuerva \cite{Gar90} (see
also \cite{GR85}) used the Maurey-Rosenthal cycle of ideas to create a
crucial link between functional analysis and harmonic analysis --
Rubio's $L_p$-boundedness principle for singular integral operators.

Continuing and improving their work, the aim of this paper is to extend the
Maurey-Rosenthal theorem within the framework of quasi K\"othe and Banach
function
spaces over measure spaces (examples are Lorentz and Orlicz spaces); we
prove a general theorem which combines the Maurey-Rosenthal theorem with
the Grothendieck-Pietsch domination theorem, with Krivine's factorization
theory for operators acting between Banach lattices, and includes
many
scalar and vector-valued variants of the original Maurey-Rosenthal theorem
as special cases.

\bsk
We shall use standard notation and notions from Banach space theory, as
presented e.g. in \cite{DJT95} or \cite{LT79}.
If $E$ is a Banach space over the scalars $\K = \R$ or $\C$, then $B_E$
denotes its closed unit ball and $E'$
its continuous dual; Banach lattices $X$ by definition are real. The term
``operator'' stands for a bounded and linear mapping between (quasi) normed
spaces. But note that most of our results can be formulated for operators
$T: E\lra F$ which are only homogeneous (in the sense that $T(\lambda x) =
\lambda T(x)$ for all $\lambda\ge 0$, $x\in E$)
and not necessarily additive -- in such a case the term ``homogeneous operator''
will be used. Recall that
for example in harmonic analysis many nonlinear operators like
square or maximal operators naturally arise.

\bn
Acknowledgement: The author thanks Marius Junge (Kiel) for various fruitful
discussions on the topic of this paper.

\section{Powers of $r$-convex quasi Banach function spaces}

In this first section we fix some terminology and lemmata on so-called
quasi K\"othe function spaces.

Let $(\Omega,\Sigma,\mu)$ be a $\sigma$-finite and complete measure space,
and denote all $\mu$-a.e. equivalence classes of real-valued measurable
functions on $\Omega$ by $L_0(\mu)$.
A quasi normed
space $\bl X(\mu),\|\cdot\|_X\br$ of functions in $L_0(\mu)$ is said to
be a quasi K\"othe function space if
it satisfies the following three conditions:
\begin{itemize}
\item[(I)]
If $|x| \le |y|$ on $\Omega$, with $x\in L_0(\mu)$ and $y\in
X(\mu)$, then
$x\in X(\mu)$  and $\|x\|_X \le \|y\|_X$.
\item[(II)]
There is some $0 < t < \infty$ such that for all
$x,y\in X(\mu)$
\[
  \bdv\bl |x|^t+|y|^t\br^{1/t}\bdv_X \le \bl\|x\|^t_X+\|y\|^t_X\br^{1/t}.
\]
\item[(III)]
The support of $X(\mu)$ (i.e., the smallest set in $\Omega$
which contains $\mu$-a.e. all supports of functions in $X(\mu)$)
equals $\Omega$, and moreover $X(\mu)$ satisfies the so-called Fatou
property:
$\|x_n\|_X \lra \|x\|_X$ for non negative functions $x_n$, $x\in
X(\mu)$ such that $x_n\uparrow x$ pointwise $\mu$-a.e.
\end{itemize}

Complete quasi K\"othe function spaces are called quasi
Banach function spaces -- if $\|\cdot\|_X$ is a (complete)
norm or, equivalently, (II) holds for $t=1$, then we shortly
speak of (Banach) K\"othe function spases.
See e.g. \cite{BS88}, \cite[II]{LT79}, \cite{Me91}
and
\cite{Schw84} for information on function spaces -- but note that in
all these references the definitions are slightly different.

Clearly, every (quasi) K\"othe function space can be interpreted as a
(quasi) normed lattice. We mention that already under weak
additional assumptions normed lattices  are order isomorphic
to K\"othe function spaces (see e.g.\ \cite{F67} and
\cite[II,1.b.14]{LT79}).

A quasi K\"othe function space $X(\mu)$ is said to be
$\sigma$-order continuous if $\|x_n\|_X \lra 0$ for each
sequence $(x_n)$ in $X(\mu)$ with $0\le x_n \downarrow 0$ $\mu$-a.e..

Define for a K\"othe function space $X(\mu)$ the intersection of the order
continuous dual and topological dual:
\[
   X^\times(\mu) := \Blk y\in L_0(\mu)  \Bmv
               \|y\|_{X^\times}:= \sup_{\|x\|_X \le 1} |\int xy\,d\mu |
               < \infty \Brk\, .
\]
The following well-known result will be important
(for $(1)$ see e.g. \cite[2.6.4 and 2.4.19]{Me91}, for $(2)$
\cite[2.4.21]{Me91} and \cite[13.5]{Schw84}, and for $(3)$
again \cite[2.6.4]{Me91}):

\begin{lemma} \label{lemma-1}
Let $X(\mu)$ be a K\"othe function space.

\begin{itemize}
\item[{\rm (1)}]
$(X^\times (\mu),\|\cdot\|_{X^\times })$ is a Banach function space.
\item[{\rm (2)}]
$\|x\|_X = \sup_{\|y\|_{X^\times } \le 1} | \int xy d \mu |\;$ for all
$x\in X(\mu)$.
\item[{\rm (3)}]
$X^\times (\mu)$=$X'(\mu)$ whenever $X(\mu)$ is $\sigma$-order
continuous.
\end{itemize}
\end{lemma}

For $0 < r < \infty$ the $r$-th power of a quasi K\"othe function space
$(X(\mu), \|\cdot\|_X)$ is defined to be
\[
   X^r(\mu) := \blk x\in L_0(\mu) \bmv |x|^{1/r} \in X(\mu)\brk\, ;
\]
together with the quasi norm
\[
   \|x\|_{X^r} := \bdv\,|x|^{1/r}\bdv^r_X\, , \q x\in X^r(\mu)
\]
this vector space obviously forms a quasi K\"othe function space
((II) then holds for $t/r$ instead of $t$). Observe that $X^r(\mu)$ is
$\sigma$-order continuous whenever $X(\mu)$ is.

We say (in analogy to the theory of Banach lattices) that a quasi
K\"othe function space $X(\mu)$ is $r$-convex $(0 < r < \infty)$ if
there is a constant $c\ge 0$ such that for all $x_1,\ldots,x_n\in X(\mu)$
\[
   \bdv\bl\tsum |x_k|^r\br^{1/r}\bdv_X \le c\bl\tsum\|x_k\|^r_X\br^{1/r}
     \, ,
\]
and $r$-concave if there is a $c\ge 0$ such that for all $x_1,\ldots,
x_n\in X(\mu)$
\[
   \bl\tsum \|x_k\|^r_X\br^{1/r} \le c \bdv\bl\tsum |x_k|^r\br^{1/r}\bdv_X
      \, ;
\]
$M^{(r)}(X)$ and $M_{(r)}(X)$ stand for the best constants, respectively
(and, as usual, we define these constants to be $\infty$ if $X$ does not
have the corresponding property).

\begin{lemma} \label{lemma-2} \text{\hfill}

\begin{itemize}
\item[{\rm (1)}]
For $0 < t,r < \infty$ and every quasi K\"othe function space
$X(\mu)$
\[
   M^{(r/t)} \bl X^t\br = M^{(r)}\bl X\br^t\,  \q\text{and}\q M_{(r/t)}
              \bl X^t\br = M_{(r)} \bl X\br^t\, .
\]
\item[{\rm (2)}]
For $1 \le r < \infty$ and every K\"othe function space
$X(\mu)$
\[
   M^{(r)} \bl X^\times\br = M_{(r')}\bl X\br\,  \q\text{and}\q M^{(r)}
\bl X^\times\br = M_{(r')}\bl X\br\, .
\]
\end{itemize}
\end{lemma}

Statement (1) is an immediate consequence of the definitions. The
proof of (2) follows from the formulas
\begin{align*}
  \bl\tsum \|x_k \|^r_X\br^{1/r}  & = \sup \bglk \bv \int \tsum
    x_k\, y_k\, d\mu\bv \Bmv \bl\tsum\|y_k\|^{r'}_{X^\times}\br^{1/r'} \le
    1\bgrk\\[1ex]
  \bdv\bl\tsum |x_k|^r\br^{1/r}\bdv_X  & =  \sup \bglk \bv\int \tsum
    x_k\, y_k\, d\mu \bv \Bmv \bdv\bl\tsum
    |y_k|^{r'}\br^{1/r'}\bdv_{X^\times} \le 1\bgrk
\end{align*}
which are consequences of the isometric embeddings
\[
\begin{array}{lclclcl}
  \e^n_r(X)  & \hra  & \e^n_r(X^{\times \times}) &
    \hra  & \e^n_r\bl(X^\times)'\br  & = &
    \bl\e^n_{r'}(X^\times)\br'\\[1ex]
  X(\e^n_r) & \hra & X^{\times\times}(\e^n_r) & \hra & (X^\times)'(\e^n_r)
     & = & \bl X^\times(\e^n_{r'})\br'
\end{array}
\]
(see Lemma~\ref{lemma-1}, (2) and \cite[II, p.~47]{LT79}, in particular
for the notation). We now clarify when the quasi norm $\|\cdot\|_{X^r}$ on
$X^r(\mu)$ is
equivalent to a norm. If $X(\mu)$ is $r$-convex, then for $x\in X^r(\mu)$
and $x_1,\ldots,x_n \in X^r(\mu)$ with $|x| \le \sum^n_{i=1} |x_i|$
\begin{align*}
   \|x\|_{X^r} = \bdv\, |x|^{1/r}\bdv^r_X
      & \le\bdv\bl \tsum \bl |x_i|^{1/r}\br^r\br^{1/r}\bdv^r_X\\[1ex]
      & \le M^{(r)}(X)^r \, \tsum \bdv\, |x_i|^{1/r}\bdv^r_X = M^{(r)}
        (X)^r \, \tsum \|x_i\|_{X^r}\, ,
\end{align*}
hence for the lattice norm
\[
   \tv x\tv_{X^r} := \inf \bglk \sum^n_{i=1} \|x_i\|_{X^r} \Bmv n\in
                      \N,\; |x| \le \sum^n_{i=1} |x_i|\bgrk\, ,
                      \q x\in X^r(\mu)
\]
the following result holds true:

\begin{lemma} \label{lemma-3}
Let $\bl X(\mu), \|\cdot\|_X\br$ be a quasi K\"othe function space
which for $0 < r < \infty$ is $r$-convex. Then $\bl X^r(\mu),\tv \cdot
\tv_{X^r}\br$ is a (normed) K\"othe function space, and on $X^r(\mu)$
\begin{equation}\label{eq-1.1}
   M^{(r)}(X)^{-r} \|\cdot\|_{X^r} \le \tv\cdot\tv_{X^r} \le
      \|\cdot\|_{X^r} \, .
\end{equation}
\end{lemma}

Note that every quasi K\"othe function space $X(\mu)$  by condition
(II) and induction is $t$-convex with constant $1$ for
some $0 < t < \infty$, hence by the preceding lemma $\bl
X^t(\mu),\|\cdot\|_{X^t}\br$ is a (normed) K\"othe function space and
\[
   \bl X(\mu),\|\cdot\|_X\br = \bl\bl X^t(\mu), \|\cdot\|_{X^t}\br^{1/t},
                                \|\cdot\|_{(X^t)^{1/t}} \br
\]
holds isometrically -- the quasi K\"othe function spaces are the powers
of (normed) K\"othe function spaces.

Recall that every Banach lattice is 1-convex  and $\infty$-concave and that
the properties ``$r$-convexity'' and ``$r$-concavity'' for $1\le r \le
\infty$ are ``decreasing and increasing  in $r$'', respectively
(\cite[II, 1.d.5]{LT79}). The argument usually given is strongly based on
duality, hence the following lemma needs an alternative proof.

\begin{lemma} \label{lemma-4}
Let $0 < t < r < \infty$. Then each $r$-convex quasi K\"othe function
space
$X(\mu)$ is $t$-convex, and
\[
   M^{(t)}(X) \le M^{(r)}(X)\, .
\]
In particular, K\"othe function spaces $X(\mu)$ being
1-convex, are $t$-convex for each $0< t < 1$.
\end{lemma}

\begin{proof}
Consider the norm $\tv\cdot\tv_{X^r}$ on $X^r(\mu)$. We will show that
on
$X^t(\mu)$
\[
   p_t(x) := \btv\,|x|^{r/t}\btv^{t/r}_{X^r}
\]
satisfies the triangle inequality. Since by \eqref{eq-1.1} on $X^t(\mu)$
\[
   M^{(r)}(X)^{-t} \|\cdot\|_{X^t} \le p_t(\cdot) \le \|\cdot\|_{X^t}\, ,
\]
the conclusion then is a simple consequence of Lemma~\ref{lemma-2}:
\begin{align*}
   M^{(t)}(X) & = M^{(t)}\bl\bl(X^t,\|\cdot\|_{X^t})^{1/t},
                  \|\cdot\|_{(X^t)^{1/t}}\br\br\\[1ex]
              & \le M^{(1)}\bl\bl X^t,\|\cdot\|_{X^t}\br\br^{1/t}
                \le M^{(r)}(X)\, .
\end{align*}
For the proof of the triangle inequality for $p_t$ show -- in a first step
and in complete analogy to the proof of the usual H\"older inequality --
that for each K\"othe (!) function space $Y(\nu)$ and $0 < u,v \le 1$ with
$u+v = 1$
\[
   Y^u\cdot Y^v \subset Y\, , \q \|xy\|_Y \le \|x\|_{Y^u}\, \|y\|_{Y^v}
\]
(see also \cite{GR85}). In particular, we obtain with $\frac{t}{r} +
\frac{r-t}{r} = 1$ that
\begin{align*}
    & (X^r)^{t/r} \cdot (X^r)^{(r-t)/r} \subset X^r\\[1ex]
    & \tv xy\tv_{X^r} \le \btv\, |x|^{r/t}\btv^{t/r}_{X^r} \; \btv\,
      |y|^{r/(r-t)}\btv^{(r-t)/r}_{X^r} \q\text{ for } x\in (X^r)^{t/r}, \;
      y\in (X^r)^{(r-t)/r}\, .
\end{align*}
Now simulate for $x,y \in X^t$ the proof of the Minkowski inequality:
\begin{align*}
   p_t(x+y)^{r/t} & = \btv\, |x+y|^{r/t} \btv _{X^r} = \btv\,
                      |x+y|\,|x+y|^{(r-t)/t}\btv_{X^r}\\[1ex]
                  & \le \btv\,|x|\,|x+y|^{(r-t)/t}\btv _{X^r}+ \btv\,
                    |y|\,|x+y|^{(r-t)/t}\btv_{X^r}\\[1ex]
                  & \le \btv\,|x|^{r/t}\btv^{t/r}_{X^r}
                    \btv\,|x+y|^{r/t}\btv^{(r-t)/r}_{X^r} + \cdots\\[1ex]
                  & = p_t(x)\, p_t(x+y)^{r/t-1} + p_t(y)\, p_t(x+y)^{r/t-1}
                     \, .
\end{align*}
This obviously completes the proof.
\end{proof}

The last lemma needed is

\begin{lemma}\label{lemma-5}
For $0 < r < \infty$ let $X(\mu)$ be an $r$-convex quasi K\"othe function
space. Then $\bl X^{r/2}(\mu), \tv\cdot\tv_{X^{r/2}}\br^\times$ is
$\sigma$-order continuous.
\end{lemma}

\begin{proof}
$X(\mu)$ by Lemma~\ref{lemma-4} is $r/2$-convex, hence $Y := \bl
X^{r/2}, \tv\cdot\tv_{X^{r/2}}\br$ by Lemma~\ref{lemma-3} is a K\"othe
function space. Moreover, $Y$ and
\[
   Z:= \bl (X^r,\tv\cdot\tv_{X^r})^{1/2},
        \|\cdot\|_{(X^r,\tv\cdot\tv_{X^r})^{1/2}}\br
\]
by \eqref{eq-1.1} are isomorphic K\"othe function spaces, and $Z$ by
Lemma~\ref{lemma-2}, (1) is 2-convex
(since the space $(X^r,\tv\cdot\tv_{X^r})$ is 1-convex). Assume now that
$Y^\times$ is not $\sigma$-order continuous -- then it contains
$\e_\infty$ as a topological subspace (note first that Banach function
spaces are $\sigma$-complete and see
e.g. \cite[II, 1.a.7]{LT79} or \cite[3.7]{Schw84}). But this contradicts
the fact that $Y^\times$ by Lemma 2, (2) is 2-concave.
\end{proof}

\bn
Our basic examples are:

(1) For $0 < p_1,p_2 \le \infty$ denote by $L_{(p_1,p_2)} (\mu_1\ot\mu_2)$
the space of all (equivalence classes of) measurable real-valued functions
$f$ on $\Omega _1 \times \Omega _2 $ such that
\[
   \|f\|_{(p_1,p_2)} := \bgl\int\bgl\int \bv f(w_1,w_2)\bv^{p_2}\,
                        d\mu_2(w_2)\bgr^{p_1/p_2} d\mu(w_1)\bgr^{1/p_1}
\]
(with the obvious modification if $p_1$ or $p_2 = \infty$). This gives a
quasi Banach function space over $\mu_1\ot \mu_2$ which is normed whenever
$1\le p_1,p_2 \le \infty$, and by the continuous triangle inequality it is
$\min(p_1,p_2)$-convex and $\max(q_1,q_2)$-concave with constants 1. For
$p_1 = p_2$ one gets the usual $L_p$'s. $L_\infty(\mu)$ being
$\infty$-convex, is $r$-convex for all $0 < r < \infty$. Obviously,
\[
   L_{(p_1,p_2)}(\mu_1\ot \mu_2)^r = L_{(p_1/r,p_2/r)}(\mu_1\ot\mu_2 )
        \q\text{ for } 0 < r < \infty, \; 0 < p_1,p_2 \le \infty\, .
\]

(2) For $0 < p_1 < \infty$ and $0 < p_2 \le \infty$ the Lorentz function
spaces $L_{p_1,p_2}(\mu)$ (see e.g. \cite{BS88} or \cite{Koen86} for the
definition) form quasi Banach function spaces over $\mu$; recall that these
spaces are normable whenever $1 < p_1 < \infty$, $1\le p_2 \le \infty$.
Again it is straightforward to check that
\[
   L_{p_1,p_2}(\mu)^r = L_{p_1/r,p_2/r}(\mu) \q\text{ for } 0 < p_1,\;
                        r < \infty, \; 0 < p_2 \le \infty\, .
\]
Convexity and concavity of Lorentz spaces for $1 < p_1 < \infty$, $0 < p_2
\le \infty$ and non-atomic $\mu$ were studied in Creekmore \cite{Cr81}. The
following
arguments handle the general case and seem to be easier: let $0 < r < p_1 <
\infty$ and $0 < r \le p_2 \le \infty$, and choose $0 < s < \infty$ such
that $1 < sr < \min (sp_1,sp_2)$. Then by Lemma~\ref{lemma-2}
\begin{align*}
   M^{(r)}(L_{p_1,p_2})  & = M^{(sr/s)}\bl L^s_{sp_1,sp_2}\br\\[1ex]
                         & \le M^{(sr)}\bl L_{sp_1,sp_2}\br^s\\[1ex]
                         & = M^{(1\bs 1/sr)}
                            \Bl L^{1/sr}_{sp_1/sr,sp_2/sr}\Br^s\\[1ex]
                         & \le M^{(1)}\bl L_{sp_1/sr,sp_2/sr}\br^{1/r}  <
                            \infty
\end{align*}
(since the latter space is normable). Similarly,
\[
   M_{(r)}(L_{p_1,p_2}) < \infty \q\text{ for } 0 < p_1 < r < \infty,\;
                                0 < p_2 \le r\, ;
\]
for $ 1 < p_1,p_2$ use the result on convexity together with duality, and
for arbitrary $0 < p_1,p_2$ Lemma~\ref{lemma-2} as above.

(3) Orlicz function spaces $L_\varphi(\mu)$ are Banach function spaces in
the above sense, and convexity and concavity can be characterized in terms
of $\varphi$ (see \cite[II]{LT79}). Clearly, for $0 < r <
\infty$
\[
   L_\varphi(\mu)^r = L_{\varphi(\sqrt[r]{\cdot})}(\mu) := \Blk f \Bmv
                      \exists \lambda > 0: \int \varphi
                      \Bl\frac{|f|^{1/r}}{\lambda^{1/r}}\Br\, d\mu <
                      \infty\Brk \, .
\]
It can easily be seen that $L_p$, $L_{(p_1,p_2)}$ and $L_{p_1,p_2}$ are
$\sigma$-order continuous whenever $0 < p, p_1, p_2 < \infty$; order
continuity of Orlicz spaces $L_\varphi$ can be characterized in terms of
so-called $\Delta_2$-conditions of $\varphi$ (see \cite[II]{LT79}).

\section{Vector-valued norm inequalities and weighted norm inequalities
for homogeneous forms}

We call a set $U$ homogeneous whenever it carries a multiplication with
positive scalars: $U\times [0,\infty[\lra U$, $(x,\lambda) \rsa \lambda x$.
If there is a homogeneous set $U$, a quasi K\"othe function space
$X(\mu)$ and a homogeneous mapping $\phi: U\lra X(\mu)$ (i.e.,
$\phi(\lambda x) = \lambda\phi(x)$ for $\lambda\ge 0$, $x\in U$), then we
say that $\phi$ {\em represents $U$ in $X(\mu)$ homogeneously}. For two
homogeneous sets $U_1, U_2$ a form $u: U_1\times U_2 \lra \K$ is said to be
homogeneous if $u(\lambda x,y) = u(x,\lambda y) = \lambda u(x,y)$ for all
$\lambda\ge 0$,
$x\in U_1$, $y\in U_2$, and a mapping $T: U_1 \lra U_2$ is homogeneous
if $T(\lambda x) = \lambda T(x)$ for $\lambda\ge 0$, $x\in U_1$.

The following result allows to transform vector-valued norm inequalities
for forms on homogeneous sets which are homogeneously representable in quasi
K\"othe function spaces, into weighted norm inequalities (and vice versa).

\begin{theo}\label{theo-1}
For $\e = 1,2$ let $0 < r_\e < \infty$ and $\frac{1}{t} = \frac{1}{r_1} +
\frac{1}{r_2}$. Let
$u: U_1 \times U_2 \lra \K$ be a homogeneous form on homogeneous sets such
that each $U_\e$ via
$\phi_\e$ can be represented homogeneously in an $r_\e$-convex quasi
K\"othe function space $X_\e(\mu_\e)$. If $u$ satisfies
\begin{equation}\label{eq-2.1}
   \bgl \sum^n_{k=1} \bv u(x_k,y_k)\bv^t\bgr^{1/t} \le \bgdv \bgl
                         \sum^n_{k=1}
                         \bv\phi_1(x_k)\bv^{r_1}\bgr^{1/r_1}\bgdv_{X_1} \;
                         \bgdv\bgl \sum^n_{k=1} \bv\phi_2(y_k)\bv^{r_2}
                         \bgr^{1/r_2}\bgdv_{X_2}
\end{equation}
for all $x_1,\ldots,x_n \in U_1$ and $y_1,\ldots,y_n\in U_2$, then there
are two positive linear functionals
\[
   \varphi_\e: X^{r_\e}(\mu_\e) \lra \R \q\text{with}\q
     \sup_{\|x\|_{X_\e}\le 1} \varphi_\e\bl |x|^r\br^{1/r} \le
      M^{(r_\e)}(X_\e)
\]
such that for all $x\in U_1$ and $y \in U_2$
\begin{equation} \label{eq-2.2}
   |u(x,y)| \le \varphi_1\bl |\phi_1(x)|^{r_1}\br^{1/r_1}
          \varphi_2\bl |\phi_2(x)|^{r_2}\br^{1/r_2}\, .
\end{equation}
If $X_\e(\mu_\e)$ is $\sigma$-order continuous, then $\varphi_{\e}$ can be
chosen to be a function in $(X_{\e}^{r_{\e}}(\mu_\e))^{\times}$.
\end{theo}

\begin{proof}
The proof is based on a standard separation argument. Define the
weak$^*$-compact and convex set
\[
   K_{\e} := \Blk \varphi\in \bl X^{r_{\e}}_{\e},
   \tv\cdot\tv_{X_{\e}^{r_{\e}}}\br' \Bmv \varphi\ge 0,\, \|\varphi\|\le 1
   \Brk \, ,
\]
and for $\ell = 1,2$ and $x_1^{(\e)},\ldots, x_n^{(\e)}\in U_{\e}$ the
affine and continuous function
\begin{align*}
   &  \phi_{x_k^{(1)},x_k^{(2)}} : K_1\times K_2\lra \K \\[1ex]
   &  \phi_{x_k^{(1)},x_k^{(2)}} (\varphi_1, \varphi_2)
                             := \sum^2_{\e=1}\, \frac{t}{r_{\e}}\,
      M_{\e}^{r_{\e}}\varphi_{\e} \bgl \sum^n_{k=1} \bv \phi_{\e}
      \bl x^{(\e)}_k\br \bv^{r_{\e}}\bgr - \sum^n_{k=1}\bv u \bl x^{(1)}_k,
      x^{(2)}_k\br\bv^t
\end{align*}
(put $M_{\e}:= M^{(r_\e)} (X_{\e})$). Note first that the set
$\ck$ of all these functions is convex: the sum of two such
functions belongs to $\ck$, and for
$\alpha\ge 0$ and $x_1^{(\e)},\ldots, x_n^{(\e)}\in U_{\e}$
\[
   \alpha \phi_{x_k^{(1)}, x_k^{(2)}}= \phi_{\alpha^{1/r_1}x_k^{(1)},
          \alpha^{1/r_2}x_k^{(2)}}.
\]
We will now show that for each $\phi_{x_k^{(1)},x_k^{(2)}}$ there is
$(\varphi_1,\varphi_2)\in K_1\times K_2$ with
\[
   \phi_{x_k^{(1)}, x_k^{(2)}}  (\varphi_1,\varphi_2) \ge 0 \, ;
\]
indeed, by the Hahn-Banach theorem and Lemma~\ref{lemma-3} there are
$\varphi_{\e}\in K_{\e}$ such that
\[
 \varphi_{\e}\Bl \sum_k \bv \phi_{\e}\bl x^{(\e)}_k\br\bv^{r_{\e}}\Br =
 \btv \sum_k \bv \phi_{\e}\bl x_k^{(\e)}\br\bv^{r_{\e}}\btv_{X^{r_{\e}}}
 \ge \frac{1}{M_{\e}^{r_{\e}}}\,
 \Bdv \Bl \sum_k \bv \phi_{\e}\bl
 x_k^{(\e)}\br\bv^{r_{\e}}\Br^{1/r_{\e}} \Bdv^{r_{\e}}_{X_{\e}} \, ,
\]
hence (recall that $a/s +b/s'\ge a^{1/s}\, b^{1/s'}$ for $s>1$ and $a,b\ge 0$)
\begin{align*}
  \phi_{x_k^{(1)},x_k^{(2)}}(\varphi_1,\varphi_2)
&\ge \prod^2_{\e=1} M_{\e}^t
  \varphi_{\e}
  \Bl \sum_k \bv \phi_{\e}\bl x^{(\e)}_k\br\bv^{r_{\e}}\Br^{t/r_{\e}} -
  \sum_k \bv u\bl x^{(1)}_k, x^{(2)}_k\br\bv^t\\[1ex]
& \ge \prod^2_{\e=1}
  \Bdv \Bl \sum_k \bv \phi_{\e}\bl
   x_k^{(\e)}\br\bv^{r_{\e}}\Br^{1/r_{\e}} \Bdv^t_{X_{\e}} \,
   - \sum_k \bv u(x_k^{(1)}, x_k^{(2)}\br\bv^t \ge 0 \, .
\end{align*}
By Ky Fan's lemma (see e.g. \cite[9.10]{DJT95}) there is
$(\varphi_1,\varphi_2)\in K_1\times K_2$ such that
\[
   \phi_{x_k^{(1)}, x_k^{(2)}}(\varphi_1,\varphi_2)\ge 0 \text{ for all }
   n \text{ and } x_1^{(\e)},\ldots, x_n^{(\e)}\in U_{\e} \,.
\]
This easily gives the conclusion: define for $x^{(1)}\in U_1, x^{(2)}\in
U_2$
\[
   a_{\e}:= M_{\e}\varphi_{\e} \bl \bv
   \phi_{\e}(x^{(\e)})\bv^{r_{\e}}\br^{1/r_{\e}} \,.
\]
Then
\begin{align*}
       \bv u\bl x^{(1)}, x^{(2)}\br\bv^t
          & = (a_1a_2)^t \bgv\, u\bgl \frac{x^{(1)}}{a_1},
              \frac{x^{(2)}}{a_2}\bgr\bgv^t\\[1ex]
          & \le (a_1a_2)^t\, \sum^2_{\e =1}\,\frac{t}{r_{\e}}\,
             M^{r_{\e}}_{\e}
            \varphi_{\e}\bgl \bgv \phi_{\e}\Bl
            \frac{x^{(\e)}}{a_{\e}}\Br\bgv^{r_{\e}}\bgr
            = (a_1a_2)^t \, ;
\end{align*}
if $a_1=0$, then $\bv u\bl n x^{(1)}, x^{(2)} \br\bv^t \le
\frac{t}{r_2}\, M^{r_2}_2\varphi_2 \bl\bv\, \phi_2(x^{(2)})\bv^{r_2}\br$
for all $n$, hence \mbox{$u(x^{(1)},x^{(2)}) = 0$.}
\end{proof}

\begin{rem}\hfill \label{rem-1} {\rm
(1)
An easy calculation shows that \eqref{eq-2.2} implies \eqref{eq-2.1} with
an additional constant $M^{(r_1)}(X_1)\, M^{(r_2)}(X_2)$.
(2)
Recall from Lemma~\ref{lemma-4} that K\"othe function spaces are $r$-convex
for each $0 < r\le 1$.
(3)
A short look at the proof shows that an analogous result holds for
forms $u: U_1 \times \ldots \times U_n \lra \K$ which are
homogeneous in each coordinate.
(4)
$\sigma$-finiteness of the measures is only needed for the last statement
of the theorem.
(5)
If $K_{\e}$ is a compact space and $\phi_{\e}: U_{\e}\lra
\e_{\infty}(K_{\e})$ has values in $C(K_{\e})$, then by restriction
$\varphi_{\e}: \e_{\infty}(K_{\e})\lra \K$ defines a Borel measure on
$K_{\e}$ (the proof shows that in this case it is even possible to obtain a
Borel probability measure).
}
\end{rem}

Together with Lemma~\ref{lemma-3} we obtain the following
interesting

\begin{cor}\label{cor-1}
For two K\"othe function spaces $X(\mu)$ and $Y(\nu)$, and $1\le r< \infty$
let
\[
   u: X^{1/r}(\mu) \times Y^{1/r'}(\nu) \lra \K
\]
be a homogeneous form such that for all $x_1,\ldots, x_n\in X^{1/r}(\mu)$,
$y_1,\ldots, y_n \in Y^{1/r'}(\nu)$
\begin{equation}\label{eq-2.3}
   \tsum \bv u(x_k,y_k)\bv \le \bdv \bl\tsum |x_k|^r\br^{1/r}\bdv_{X^{1/r}}
   \bdv \bl \tsum | y_k|^{r'}\br^{1/r'}\bdv_{Y^{1/r'}}\, .
\end{equation}
Then there exist positive functionals
\begin{align*}
  \varphi\in X'(\mu)  & \q\text{with}\q \|\varphi\|_{X'}\le
    1\\[1ex]
  \psi\in Y'(\nu)     & \q\text{with}\q \|\psi \|_{Y'}\le 1
\end{align*}
such that for all $x\in X^{1/r}(\mu), \, y\in Y^{1/r'}(\nu)$
\[
   \bv u(x,y)\bv \le \varphi\bl |x|^r\br^{1/r} \psi \bl |y|^{r'}
   \br^{1/r'} \, .
\]
If $X(\mu)$ and $Y(\nu)$, respectively, are $\sigma$-order continuous, then
$\varphi$ and $\psi$ are functions in $X^{\times}(\mu)$ and $Y^{\times}(\nu)$,
respectively.
\end{cor}

\begin{proof}
Take $U_1 = X^{1/r}$, $U_2 = Y^{1/r'}$ and $\phi_1, \phi_2$ the identities,
and recall from Lemma~\ref{lemma-2}, (1) that $X^{1/r}$ is $r$-convex and
$Y^{1/r'}$ $r'$-convex with constants 1.
\end{proof}

For $r = 2$ every continuous and bilinear $u: X^{1/2} \times Y^{1/2}\lra
\K$ satisfies
\eqref{eq-2.1} -- this is an easy consequence of the Grothendieck-Krivine
theorem from \cite{Kr73} (see also \cite[II, 1.f.14]{LT79}).
Moreover, it can be shown that each positive continuous bilinear form
$u: X^{1/r} \times X^{1/r'}\lra \K$, i.e. $u(x,y)\ge 0$ for $x,y\ge
0$, satisfies \eqref{eq-2.3} (use \cite[II, 1.d.9]{LT79}).

\section{Vector-valued norm inequalities and weighted norm inequalities
for homogeneous operators}

The following theorem is our main result -- it is a sort of reformulation
of Theorem~\ref{theo-1} for homogeneous operators instead of forms. Before
we formulate it,
let us sketch how to reprove (part of) the original Maurey-Rosenthal
theorem on the basis of Theorem~\ref{theo-1} in order to motivate
the proof of our extension which at first glance may look strange:
take an operator $T: E\lra L_p(\nu)$ ($E$ a quasi normed space and $1
\le p < \infty$) which for
$1 < r < \infty$ satisfies for all $x_1,\ldots, x_n\in E$
\[
   \bdv \bl \tsum |Tx_k|^r \br^{1/r}\bdv_{L_p(\nu)} \le \bl \tsum
                  \|x_k\|_E^r \br^{1/r}\, .
\]
Then
$T$ defines a bilinear form
\[
   u: E\times L_p(\nu)^{\times}\lra \K, \, u(x,y) := \int Tx\cdot y\, d\nu
       \, .
\]
If we consider the ``trivial representations''
\begin{align*}
 &   \phi_1: E\lra L_r(\mu), \, \phi_1x := \| x\|_E 1\q
            \q\text{($\mu$ some Dirac measure)}\\
 &   \phi_2:  L_p(\nu)^{\times} \lra L_p(\nu)^{\times},\, \phi_2x := x\, ,
\end{align*}
then, by H\"older's inequality, $u$ satisfies \eqref{eq-2.1}. Hence, if
$1\le p< r<\infty$, then $L_p(\nu)^{\times}$ is $r'$-convex,
and Theorem~\ref{theo-1} gives a functional
\[
   0 \le \varphi\in \bl (L_p(\nu)^{\times})^{r'} \br' \q\text{ with
        }\q\|\varphi\|\le 1
\]
such that for all $x\in E, \, y\in L_p(\nu)^{\times}$
\[
   \bv u(x,y)\bv \le \|x\|_E \, \varphi\bl |y|^{r'}\br^{1/r'} \, .
\]
For $p> 1$ this functional $\varphi$ is a function in
$L_{(p'/{r'})'}(\nu)$,
and it is not difficult to show (for details see the next proof) that with
$\omega := \varphi^{r/r'}$ for all $x\in E$
\[
   \int \frac{|Tx|^r}{\omega} d\nu\le \|x\|_E \, .
\]
Hence, for $1< p < r < \infty$ we have reproved the
original Maurey-Rosenthal theorem; for $p=1$
we essentially do the same -- but we need some scaling trick in order to
get a {\it function} $\varphi\in L_1(\nu)$ instead of only a measure in
$\bl (L_1(\nu)^{\times})^{r'}\br' = L_{\infty}(\nu)'$. This trick will even
cover the general case $0 < p < r < \infty$ although the above argument is
strongly based on duality.

\begin{theo} \label{theo-2}
Let $T: U\lra V$ be a homogeneous operator where $U$ and $V$ are vector
spaces (or only
homogeneous sets) which via $\phi$ and $\psi$ are represented homogeneously
in quasi
K\"othe function spaces $X(\mu)$ and $Y(\nu)$, respectively. For $0 < r <
\infty$ let $X(\mu)$ be $r$-convex, and $Y(\nu)$ $r$-concave.
If $T$ satisfies
\begin{equation}\label{eq-3.1}
  \bdv \bl \tsum \bv \psi(Tx_k)\bv^r\br^{1/r}\bdv_Y \le
  \bdv \bl \tsum \bv \phi (x_k)\bv^r\br^{1/r}\bdv_X
\end{equation}
for all $x_1,\ldots,x_n \in U$, then there are
\begin{align}\notag
    0\le \varphi: X^r\lra \R & \q\text{linear with}\q
     \sup_{\| x \|_{X}\le 1}\,  \varphi (|x|^r)^{1/r}
     \le M^{(r)}(X)\\[1ex]
 \label{eq-3.2}
    0\le \omega_2\in L_0(\nu)\hspace*{3.4mm} &
    \q\text{with\hspace*{1.1cm}}\q
    \sup_{\| y\|_{L_r(\nu)}\le 1}\,
    \bdv \omega_2^{1/r}y\|_{Y}\le M_{(r)}(Y)
\end{align}
such that for all $x\in U$
\begin{equation} \label{eq-3.3}
  \int \frac{|\psi(Tx)|^r}{\omega_2}\, d\nu \le \varphi
   \bl |\phi(x)|^r\br \, .
\end{equation}
If $X(\mu)$ is $\sigma$-order continuous, then
there is even a function
\begin{equation}\label{eq-3.4}
     0\le \omega_1\in L_0(\mu) \q \text{with} \q \sup_{\| x\|_X\le 1} \bdv
     \omega_1^{1/r}x\bdv_{L_r(\mu)}\le M^{(r)}(X)
\end{equation}
such that for all $x\in U$
\[
   \int \frac{\bv \psi(Tx)\bv^r}{\omega_2}\, d\nu \le \int \bv \phi
   (x)\bv^r \omega_1\, d\mu\, .
\]
\end{theo}

Recall that K\"othe function spaces by Lemma~\ref{lemma-4} are $r$-convex
for all $0 < r\le 1$.

\begin{proof}
We may assume without loss of generality that $V = Y$ and $\psi$ is the
identity (otherwise replace $T$ by $\psi \ci T$). By assumption
(property (II)) $Y(\nu)$ is
$t$-convex with constant 1 for some $0 < t < \infty$, and by
Lemma~\ref{lemma-4} we may
assume that $t < r$. Consider for $s = t/2$ a new scalar multiplication on
$U$: for $\lambda\ge 0$ and $x\in U$ define
\[
   \lambda \ci x := \lambda^{1/s}x\, ;
\]
then it is easy to check that the mappings
\begin{align*}
     \phi^s & : U\lra X^s, & \hspace*{-3cm}  \phi^s(x) & :=
       |\phi(x)|^s\\
     T^s    & : U\lra Y^s,  & \hspace*{-3cm} T^s x     & :=
       |Tx|^s
\end{align*}
with respect to this multiplication
are homogeneous, and for $x_1,\ldots, x_n\in U$
\begin{equation}\label{eq-3.5}
   \bdv \bl \tsum \bv T^s x_k\bv^{r/s}\Br^{s/r}\bdv_{Y^s} \le
   \bdv \bl \tsum \bv \phi^s(x_k)\bv^{r/s}\Br^{s/r}\bdv_{X^s}\, .
\end{equation}
From now on $Y^s$ will always be endowed with its natural quasi norm
$\|\cdot\|_{Y^s}$ which by Lemma~\ref{lemma-3} is even a norm since
$M^{(s)}(Y) \le M^{(t)}(Y) \le 1$ by Lemma~\ref{lemma-4}. Hence, the
homogeneous form
\[
   u^s: U\times(Y^s)^{\times} \lra \R, \q u^s(x,y) := \int T^sx
   \cdot y \,  d\nu
\]
is defined, and
by H\"older's inequality, the inequality  \eqref{eq-3.5} on $T^s$, and the
definition of the norm on $(Y^s)^{\times}$ it satisfies
\begin{align*}
   \tsum\bv u^s(x_k,y_k)\bv
           & \le \int \tsum\bv T^s x_k\bv\, |y_k|\, d\nu\\[1ex]
           & \le \int \bl\tsum \bv T^s x_k\bv^{r/s}\br^{s/r} \bl\tsum
                 |y_k|^{(r/s)'}\br^{1/(r/s)'}\, d\nu\\[1ex]
           & \le \bdv\bl\tsum \bv T^s x_k\bv^{r/s}\br^{s/r}\bdv_{Y^s}\,
                  \bdv\bl\tsum |y_k|^{(r/s)'}
                  \br^{1/(r/s)'}\bdv_{(Y^s)^{\times}}
\end{align*}
for all $x_1,\ldots,x_n\in U$ and $y_1,\ldots,y_n \in (Y^s)^{\times}$. By
Lemma~\ref{lemma-2}
\begin{align*}
   M^{(r/s)}(X^s) & = M^{(r)}(X)^s < \infty\\[1ex]
   M^{((r/s)')}\bl (Y^s)^\times\br  & = M_{(r/s)}(Y^s) = M_{(r)}(Y)^s <
                                        \infty \, ,
\end{align*}
hence it follows from Theorem~\ref{theo-1} that there are two functionals
\begin{align*}
   0\le \varphi_1: \bl X^s,\|\cdot\|_{X^s}\br^{r/s} \lra \R\,,  &
         \q\sup_{\|x\|_{X^s}\le 1} \bv \varphi_1(|x|^{r/s})\bv^{s/r}
             \le M^{(r)}(X)^s\\[1ex]
   0\le \varphi_2: \bl(Y^s)^\times, \|\cdot\|_{(Y^s)^\times}\br^{(r/s)'}
             \lra  \R\,, &
         \q\sup_{\|y\|_{(Y^s)^\times}\le 1}
             \bv\varphi_2(|y|^{(r/s)'})\bv^{1/(r/s)'}
             \le M_{(r)}(Y)^s
\end{align*}
such that for all $x\in U$, $y\in (Y^s)^\times$
\begin{equation} \label{eq-3.6}
   |u^s(x,y)| \le \varphi_1\bl |\phi^s(x)|^{r/s}\br^{s/r}\, \varphi_2\bl
                   |y|^{(r/s)'}\br^{1/(r/s)'}\, .
\end{equation}
Clearly, for $\varphi:= M^{(r)}(X)^{-r}\varphi_1$
\begin{equation} \label{eq-3.7}
   0\le \varphi: X \lra \R\, , \q \sup_{\|x\|_X\le 1} \varphi(|x|^r)^{1/r}
                 \le 1\, .
\end{equation}
On the other hand (Lemma~\ref{lemma-3})
\[
   \varphi_2 \in \blb\bl(Y^s)^\times, \|\cdot\|_{(Y^s)^\times}\br^{(r/s)'},
    \tv\cdot\tv\ldots\brb' \, ,
\]
but since $(Y^s)^\times$ and hence also its $(r/s)'$-th power by
Lemma~\ref{lemma-5} are $\sigma$-order continuous, we obtain from
Lemma~\ref{lemma-1}, (3)
\begin{equation} \label{eq-3.8}
   g := M_{(r)}(Y)^{-s(r/s)'} \varphi_2 \in L_0(\nu), \q
    \sup_{\|y\|_{(Y^s)^\times}\le 1} \int g|y|^{(r/s)'}\, d\nu \le 1\, .
\end{equation}
Moreover, by \eqref{eq-3.6} for all $x\in U$, $y \in (Y^s)^\times$
\[
  \int |Tx|^s y\, d\nu \le M^{(r)}(X)^{s}M_{(r)}(Y)^{s} \varphi\bl |\phi
     (x)|^r\br^{s/r}\,
     \bgl \int g|y|^{(r/s)'}\, d\nu\bgr^{1/(r/s)'} \, .
\]
It remains to give this inequality the form of \eqref{eq-3.3}:
note first that for $h := g^{1/(r/s)'}$ the multiplication operator
\[
   M_h : (Y^s)^\times \lra L_{(r/s)'}\bl \nu_{\bv [h>0]}\br
\]
by \eqref{eq-3.8} is defined, and for all $x\in U$, $y\in (Y^s)^\times$
\[
   \int\frac{|Tx|^s}{h}\, M_h y \, d\nu \le M^{(r)}(X)^{s} M_{(r)}(Y)^{s}
           \varphi\bl|\phi(x)|^r\br^{s/r}\, \|M_h y\|_{L_{(r/s)'}}
\]
(obviously, $[h=0] \subset [Tx = 0]$ $\nu$-a.e. for all $x\in U$). Since
$M_h$ has dense range (use e.g.\ the Hahn-Banach theorem), for all $x\in U$
\[
   \int\frac{|Tx|^r}{h^{r/s}}\, d\nu \le M^{(r)}(X)^{r} M_{(r)}(Y)^{r}
         \varphi\bl|\phi(x)|^r\br\, .
\]
Finally, it remains to prove the norm estimates \eqref{eq-3.2} and
\eqref{eq-3.3} for the weights
\begin{align*}
   \omega_2 &  := M^{(r)} (X)^r h^{r/s}\in L_0(\nu)\\
   \omega_1 &  := M^{(r)}(X)^r \varphi\in L_0(\mu)\q \text{ (provided } X
                   \text{ is $\sigma$-order continuous})\, .
\end{align*}
The estimate in \eqref{eq-3.8} and the fact that $B_{(Y^s)^{\times}}$ is
norming in $Y^s$ (see Lemma~\ref{lemma-1}, (2)) give
\begin{align*}
  \sup_{\| y\|_{L_r}\le 1}\|h^{1/s}y\|_Y
   & = \sup_{\| y\|_{L_{r/s}}\le 1}\, \bdv h^{1/s}|y|^{1/s}\bdv_Y\\[1ex]
   & = \sup_{\| y\|_{L_{r/s}}\le 1} \,\bdv hy\bdv^{1/s}_{Y^s}\\[1ex]
   & = \sup_{\| y\|_{L_{r/s}}\le 1} \; \sup_{\|f\|_{(Y^s)^\times}\le 1}
       \bv\int hyf\, d\nu\bv^{1/s}\\[1ex]
   & = \sup_{\| f\|_{(Y^s)^{\times}}\le 1} \,\bdv
        hf\bdv^{1/s}_{L_{(r/s)'}}
     = \sup_{\| f\|_{(Y^s)^{\times}}\le 1} \,\bdv g |f|^{(r/s)'}
        \bdv^{(1/s)(1/(r/s)'}_{L_1} \le 1\, ,
\end{align*}
and if $X$ is $\sigma$-order continuous, then \eqref{eq-3.7} gives
\[
  \sup_{\|x\|_X\le 1}\bdv \varphi^{1/r}x\bdv_{L_r}
   = \sup_{\| x\|_X\le 1} \varphi \bl |x|^r\br^{1/r} \le 1\, .
\]
This completes the proof.
\end{proof}

\begin{rem}\label{rem-2} {\rm
All assumptions of the preceding result are necessary in a sense:
(1)
A straightforward computation shows that \eqref{eq-3.3}
implies \eqref{eq-3.1} with the additional constant $M_{(r)}(Y) M^{(r)}(X)$.
(2)
In which sense the convexity assumption on $X$ and concavity assumption
on $Y$ are necessary will be seen in 4.2.
(3)
That the $\sigma$-order continuity of $X$ is unavoidable in order to obtain
a function $\omega_1$, will be shown in Remark~\ref{rem-4} in 4.3.
(4)
Clearly, a remark analogous to Remark~\ref{rem-1}, (5)
is possible.
}
\end{rem}

\section{Variants of the Maurey-Rosenthal theorem for quasi Banach function
spaces}

In this section we want to illustrate that the preceding theorem
starts living if one looks at the following natural representations
$\phi: U\lra X(\mu)$ (for simplicity we will always consider complete
spaces):

\begin{itemize}
\item[(A)]
$U=E$ a Banach space, $X = \e_{\infty}(B_{E'})$ and
\[
   \phi: E\lra \e_{\infty}(B_{E'}), \, (\phi x)(x') := x'(x)\, .
\]
\item[(B)]
$U = E$ a quasi Banach space, $X = L_r(\mu)$ (with $0 < r\le \infty$,
$\mu$ a Dirac measure) and
\[
  \phi: E\lra L_r(\mu), \q \phi x:= \| x\|_E 1 \, .
\]
\item[(C)]
$U = X$ a quasi Banach function space and the identity
\[
   \phi : X\lra X, \q \phi x:= x \, .
\]
\item[(D)]
$U = L_r(\mu,E)$ a space of Bochner integrable functions with values in a
quasi Banach space $E$, $X = L_r(\mu)$ and
\[
   \phi: L_r(\mu,E)\lra L_r(\mu), \q \phi x := \| x\|_E(\cdot)\, .
\]
\item[(E)]
More generally, take a quasi Banach function space $X(\mu)$ and a quasi
Banach space $E$. Then the vector-valued quasi Banach function space
\begin{align*}
    & X(\mu,E) := \blk x:\Omega \lra E \mid x \text{
                $\mu$-measurable, $\|x\|_E(\cdot) \in X(\mu)$}\brk\\[1ex]
    & \|x\|_{X(\mu,E)} := \bdv\, \|x\|_E(\cdot) \bdv_X
\end{align*}
is a quasi Banach space, and
\[
   \phi : X(\mu,E) \lra X(\mu), \; \phi x := \|x\|_E(\cdot)
\]
represents $X(\mu,E)$ in $X(\mu)$ homogeneously.
\end{itemize}


\subsection{The Grothendieck-Pietsch domination theorem revisited}

It can easily be seen that the fundamental Grothendieck-Pietsch
domination/factorization theorem for summing operators is a
straightforward consequence of Theorem~\ref{theo-2}. Let $T: E\lra F$
be an $r$-summing ($1\le r < \infty$) operator between Banach spaces, i.e.,
there is some constant $c\ge 0$ such that for all $x_1,\ldots, x_n\in E$
\[
   \bl \tsum \|Tx_k \|_F^r\br^{1/r} \le c \, \sup_{\| x'\|_{E'}\le 1}\, \bl
   \tsum \bv x'(x_k)\bv^r\br^{1/r}\, ;
\]
as usual $\pi_r(T) := \rinf c$ denotes the $r$-summing norm of $T$. With
the representations
\begin{align*}
    & \phi: E\lra \e_{\infty}(B_{E')}  \q\text{ as in (A)}\\
    & \psi: F\lra L_r(\nu)  \qq\text{ as in (B)}
\end{align*}
this inequality transfers into
\[
   \bdv \bl \tsum \bv \psi(Tx_k)\bv^r\br^{1/r}\bdv_{L_r(\nu)} \le
     \pi_r(T) \bdv \bl \tsum\bv
     \phi(x_k)\bv^r\br^{1/r}\bdv_{\e_{\infty}(B_{E'})}\, .
\]
Since $\e_{\infty}(B_{E'})$ is $r$-convex and $L_r(\nu)$ $r$-concave (with
constants 1), Theorem~\ref{theo-2} (see also Remark~\ref{rem-1}, (4)) gives
a positive $\varphi\in \e_{\infty}(B_{E'})'$ of norm $\le 1$ such that
\[
   \|Tx \|^r \le \pi_r(T)^r \varphi\bl \bv\la \cdot,x\ra\bv^r\br\, ;
\]
this reproves the Grothendieck-Pietsch domination theorem (clearly,
$\varphi$ defines a Borel measure on the weak$^*$-compact set $B_{E'}$).


\subsection{The Maurey-Rosenthal theorem for $r$-convex
homogeneous operators with values in $r$-concave quasi Banach function
spaces}

Let $T: E\lra Y(\nu)$ be a homogeneous operator, where $E$
is a quasi Banach space and $Y(\nu)$ a quasi Banach function space. As in
the linear and normed case, we call $T$ $r$-convex $(0 < r < \infty)$ if
there is a
$c\ge 0$ such that for all $x_1,\ldots, x_n\in E$
\[
   \bdv \bl\tsum |Tx_k|^r\br^{1/r}\bdv_Y \le c\,\bl\tsum  \|
   x_k\|^r_E\br^{1/r}\, ;
\]
the best constant $c$ is denoted by $M^{(r)}(T)$.

\bsk
A theorem of Krivine \cite{Kr73} (see \cite[II,1.d.12]{LT79}) states that
for $1\le r < \infty$ every
$r$-convex operator $T$ from a Banach space $E$ into an $r$-concave Banach
lattice $Y$ allows a factorization

\begin{center}
\hspace*{-2cm}
\special{em:linewidth 0.4pt}
\unitlength 0.60mm
\linethickness{0.4pt}
\begin{picture}(47.66,38.00)
\put(2.67,34.00){\makebox(0,0)[cc]{$E$}}
\put(42.67,34.00){\makebox(0,0)[cc]{$Y$}}
\put(42.67,4.00){\makebox(0,0)[cc]{$L_r(\nu)$}}
\put(38.00,34.33){\vector(1,0){0.2}}
\emline{7.67}{34.33}{1}{38.00}{34.33}{2}
\put(39.00,9.33){\vector(3,-2){0.2}}
\emline{4.67}{30.33}{3}{39.00}{9.33}{4}
\put(22.67,38.00){\makebox(0,0)[cc]{\footnotesize$T$}}
\put(17.33,18.00){\makebox(0,0)[cc]{\footnotesize$R$}}
\put(47.66,22.33){\makebox(0,0)[cc]{\footnotesize$S$}}
\put(42.67,30.00){\vector(0,1){0.2}}
\emline{42.67}{9.67}{5}{42.67}{30.00}{6}
\end{picture}
   \hspace*{1.4cm}\raisebox{1.6cm}{
   \parbox[t]{4cm}{
            $\nu$ some measure \\[0.3ex]
            $R,S$ operators\\[0.3ex]
            $S$ positive. } }
   \end{center}

\noi For a smaller class of $Y$'s
the Maurey-Rosenthal theorem shows that $S$ can be chosen to be better:
every $r$-convex operator $T: E\lra L_p(\nu)$ $(1\le p< r<\infty)$ has a
factorization $T = SR$ as above, where $S: L_r(\nu)\lra Y$ is a
positive multiplication operator.

\bsk
In view of the fact that every Banach lattice under weak
additional assumptions is
isomorphic to a Banach function space (see Section~1) the following result
in a sense combines these two theorems.

\begin{cor}\label{cor-2}
For $0 < r < \infty$ let $T$ be an $r$-convex homogeneous operator from a
quasi Banach space $E$ into an $r$-concave quasi Banach function space
$Y(\nu)$. Then there is a weight
\[
   0 \le \omega \in L_0(\nu) \q \text{with}\q \sup_{\| y \|_{L_r(\nu)}\le
                1}\, \bdv \omega^{1/r}y\bdv_Y \le M^{(r)}(T)\,  M_{(r)}(Y)
\]
such that for all $x\in E$
\[
   \bgl\int \frac{|Tx|^r}{\omega} \, d\nu \bgr^{1/r} \le  \| x\|_E \, .
\]
If $T$ is moreover linear and $1\le r < \infty$, then it factorizes as
follows:
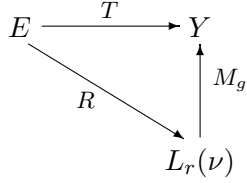
\begin{center}
\hspace*{1cm}
\special{em:linewidth 0.4pt}
\unitlength 0.60mm
\linethickness{0.4pt}
\begin{picture}(49.33,38.00)
\put(2.67,34.00){\makebox(0,0)[cc]{$E$}}
\put(42.67,34.00){\makebox(0,0)[cc]{$Y$}}
\put(42.67,4.00){\makebox(0,0)[cc]{$L_r(\nu)$}}
\put(38.00,34.33){\vector(1,0){0.2}}
\emline{7.67}{34.33}{1}{38.00}{34.33}{2}
\put(39.00,9.33){\vector(3,-2){0.2}}
\emline{4.67}{30.33}{3}{39.00}{9.33}{4}
\put(22.67,38.00){\makebox(0,0)[cc]{\footnotesize$T$}}
\put(17.33,18.00){\makebox(0,0)[cc]{\footnotesize$R$}}
\put(49.33,22.33){\makebox(0,0)[cc]{\footnotesize$M_g$}}
\put(42.67,30.00){\vector(0,1){0.2}}
\emline{42.67}{9.67}{5}{42.67}{30.00}{6}
\end{picture}
   \hspace*{1.5cm}\raisebox{1.6cm}{
   \parbox[t]{7cm}{
            $M_g$ a positive multiplication operator \\[0.3ex]
            $R$ an operator\\[0.3ex]
            $\| M_g\| \, \| R \| \le M^{(r)}(T)\, M_{(r)}(Y)$. } }
\end{center}
\end{cor}

Under the additional assumptions that $T$ is sublinear, $1<r$ and $(Y^{\times})^{r'}$
is reflexive, this result was proved by Garcia-Cuerva in
\cite[Theorem~2.9]{Gar90}.

\begin{proof}
Clearly, the proof of the first statement is an immediate consequence of
Theorem~\ref{theo-2} (represent $E$ as in (B) and $Y$ as in (C)), and
for the factorization define
\[
    g:= \omega^{1/r}\q \text{and}\q Rx := \frac{Tx}{g},\,
    x\in E\, .
\]
\end{proof}

Corollary~\ref{cor-1} is of special interest if $E$ is a Banach space
such that for {\em all\/} $r$-concave Banach function spaces $Y$ {\em every\/}
operator $T:E\lra Y$ is $r$-convex.
In the following corollary it is combined with several results and
techniques from \cite{Kr73} and \cite{Mau73,Mau74}.

\begin{cor}\label{cor-3}
Let $1\le r<\infty$ and $E$ be a Banach space. Then the following are
equivalent:
\begin{itemize}
\item[{\rm (1)}]
Every $r$-summing operator on $E'$ is 1-summing (=: the identity on $E'$
is $(r,1)$-mixing).
\item[{\rm (2)}]
Every operator $T: E\lra Y$, $Y$ a Banach lattice, is $r$-convex.
\item[{\rm (3)}]
For every operator $T: E\lra Y(\nu)$, $Y(\nu)$ an $r$-concave Banach
function space, there is some $0\le \omega\in L_0(\nu)$
with $\sup_{\| y \|_{L_r(\nu)}\le 1}\, \| \omega^{1/r}y\|_Y < \infty$
such that for all $x\in E$
\[
   \Bl \int \frac{| Tx|^r}{\omega}\, d\nu\Br^{1/r}\le \| x \|_E\, .
\]
\end{itemize}
\end{cor}

\noi
We only sketch the

\begin{proof}
If one replaces the $Y$'s in (3) by the class of all $L_p$'s
$(1\le p<r< \infty)$, then the
equivalence $(1) \carle (3)$ is a well-known result of \cite{Mau74} (see
also
\cite[32.6 and 32.5]{DF93}). Since $(2) \car (3)$ follows from
Corollary~\ref{cor-2},
it remains to check the implication $(1) \car (2)$: note first that by
Maurey's result
(just mentioned) there is a $c\ge 0$ such that
$ M^{(r)}(T) \le c\, \|T\| $ for all operators $T: E\lra \e_1$.
The use of conditional expectation operators shows that this
inequality also holds for all $T: E\lra L_1(\nu)$, $\nu$ an arbitrary
measure.
But then Krivine's localization technique from \cite{Kr73} (see also
\cite[II, the proof of 1.f.14]{LT79}) assures that every
$T: E\lra Y$, $Y$ an abitrary Banach lattice, is $r$-convex.
\end{proof}

\begin{rem}\label{rem-3} {\rm
Here we collect some well-known statements on the class of all Banach
spaces
$E$ satisfying (1):
\begin{itemize}
\item[(a)]
Each $E$ such that $E'$ has cotype 2
for $r=2$ satisfies (1). Conversely, if $E$ for $r=2$ satisfies (1), then
$E'$ has cotype $2+\varepsilon$ for all $\varepsilon> 0$.
Moreover, for $1<q<2$ each $E$ such that $E'$ has cotype $q'$ satisfies
(1) for $r = q + \varepsilon$, $\varepsilon > 0$.
See \cite{Mau74} and \cite{MP76}.
\item[(b)]
For $1\le r<2$ each $E$ with stable type $r$ satisfies (1) (see again
\cite{Mau74}). Recall that $L_q(\mu)$ for $1\le r< q\le 2$ has stable type
$r$.
\item[(c)]
Each Banach lattice $E$ with (1) is $r$-convex: by trace duality (1) holds
if and only if every $T: \e_{\infty}\lra E'$ is
$r$-summing (see e.g. \cite[32.2]{DF93})
which by a result of Maurey from \cite{Mau73} (see also
\cite[II, 1.d.10]{LT79})
implies that $E'$ is $r'$-concave, hence $E$ is $r$-convex.
For $r=2$ a Banach lattice
satisfies (1) if and only if $E'$ has cotype 2 if and only if $E$ is
2-convex (see \cite{Mau73} and \cite[II, 1.f.16]{LT79}).
\end{itemize} }
\end{rem}

\bsk
Another important result of \cite{Kr73} (see also \cite[II,1.d.9]{LT79})
shows
that for all positive operators $T: E\lra Y$ between two Banach lattices
$E$ and $Y$
\[
  \bdv \bl \tsum \bv Tx_k\bv^r\br^{1/r}\bdv_Y \le
   \| T \| \bdv \bl \tsum \bv  x_k\bv^r\br^{1/r}\bdv_E
  \text{ for all } x_1,\ldots, x_n\in E \, .
\]
Hence, if $E$ is $r$-convex, then $M^{(r)}(T)\le M^{(r)}(E) \| T \|$ which
together with Corollary~\ref{cor-2} proves the implications (1) $\car$
(2) $\car$
(3) of the following result -- the remaining implication (3) $\car$ (1)
will be shown in a more general context in \cite{DJII}.

\begin{cor}\label{cor-4}
Let $1\le r<\infty$ and $E$ be a Banach lattice. Then the following are
equivalent:
\begin{itemize}
\item[{\rm (1)}]
$E$ is $r$-convex.
\item[{\rm (2)}]
Every positive operator $T: E\lra Y$, $Y$ a Banach lattice, is
$r$-convex.
\item[{\rm (3)}]
Each positive operator $T: E\lra Y(\nu)$, $Y(\nu)$ an $r$-concave Banach
function space, satisfies {\rm (3)} of Corollary~{\rm
\ref{cor-2}}.
\end{itemize}
\end{cor}


\subsection{The Maurey-Rosenthal theorem for $r$-concave homogeneous
operators on $r$-convex quasi Banach function spaces}

Let us now consider the dual result of Corollary~\ref{cor-1}. A
homogeneous operator
$T: X(\mu)\lra F$, $X(\mu)$ a quasi Banach function space and $F$ a
quasi Banach space, is said
to be $r$-concave $(0 < r < \infty)$ if there is a $c\ge 0$ such that for
all $x_1,\ldots, x_n\in X$
\[
   \bl \tsum \| Tx_k\|_F^r\br^{1/r}\le c\,\bdv \bl \tsum |x_k|^r
       \br^{1/r}\bdv_X\, ,
\]
and the best $c$ in this inequality is denoted by $M_{(r)}(T)$.

\begin{cor}\label{cor-5}
For $0 <r< \infty$ let $T$ be an $r$-concave homogeneous operator from
an $r$-convex quasi Banach function space $X(\mu)$  into
a quasi Banach space $F$. Then there is a linear functional
$0\le \varphi: X^r(\mu)\lra \R$ with $\sup_{\| x\|_{X}\le 1}
\varphi \bl |x|^r\br^{1/r}\le  M_{(r)}(T) M^{(r)}(X)$
and such that for all $x\in X(\mu)$
\[
   \|Tx\|_F \le \varphi(|x|^r)^{1/r}.
\]
If $X$ is moreover $\sigma$-order continuous, $T$ linear and $1\le r <
\infty$, then $T$ factorizes as follows:

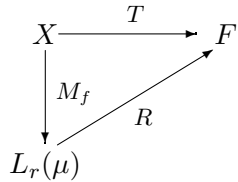
\begin{center}
\hspace*{0.5cm}
\special{em:linewidth 0.4pt}
\unitlength 0.60mm
\linethickness{0.4pt}
\begin{picture}(51.00,36.00)
\put(51.00,33.33){\makebox(0,0)[cc]{$F$}}
\put(11.00,33.33){\makebox(0,0)[cc]{$X$}}
\put(11.00,4.00){\makebox(0,0)[cc]{$L_r(\mu)$}}
\put(31.00,36.00){\makebox(0,0)[cb]{\footnotesize$T$}}
\put(32.67,16.67){\makebox(0,0)[ct]{\footnotesize$R$}}
\put(13.34,20.00){\makebox(0,0)[lc]{\footnotesize$M_f$}}
\put(11.00,9.34){\vector(0,-1){0.2}}
\emline{11.00}{29.67}{1}{11.00}{9.34}{2}
\put(48.00,29.67){\vector(3,2){0.2}}
\emline{13.67}{8.67}{3}{48.00}{29.67}{4}
\put(44.33,33.34){\vector(1,0){0.2}}
\emline{14.00}{33.34}{5}{44.33}{33.34}{6}
\end{picture}
   \hspace*{1cm}\raisebox{1.6cm}{
   \parbox[t]{7cm}{$M_f$ a positive multiplication operator\\[0.6ex]
            $R$ an operator\\[0.6ex]
            $\| M_f\|\, \| R\| \le M_{(r)}(T) M^{(r)}(X)$ .} }
   \end{center}
\end{cor}

\begin{proof}
Again, represent $X$ as in (C), $F$ as in (B) and use Theorem~\ref{theo-2}.
It remains to prove that for $\sigma$-order continuous $X$ linear $T$
factorize as indicated: with the function $\omega$ representing $\varphi\in
(X^r)'= (X^r)^{\times}$ define $f:= \omega^{1/r}$, and
\[
  R_0: \bl \ran \, M_f, \, \|\cdot \|_{L_r}\br \lra F, \, R_o(fx) := Tx
       \, .
\]
Then $\| M_f\| \le M_{(r)}(T)\, M^{(r)}(X)$, and
since $\ov{\ran}\, M_f = \blk x\in L_r(\mu) \bmv x_{\bv [f=0]} = 0\brk$ and
$L_r(\mu) = \ov{\ran}\, M_f \oplus_r \blk x\bmv x_{\bv [f>0]} = 0\brk$, the
operator $R_0$ has an extension $R$ to all of $L_r(\mu)$ with $\| R\| \le
1$. Clearly, $T= R M_f$.
\end{proof}

\begin{rem}\label{rem-4} {\rm
For $1\le r < \infty$ and sublinear $T$ between Banach spaces this result
is due to Garcia-Cuerva \cite[Theorem~2.3]{Gar90}; for
$r=2$ and linear $T$ the second part is mentioned
without proof in \cite{Ru87} -- but note
that there the assumption on the $\sigma$-order continuity of $X$ is
missing which makes the statement false:
assume that each continuous functional $\varphi: X\lra \K$ factorizes as
above: $\varphi = RM_f$. Represent $R$ by a function $g\in L_{r'}(\mu)$.
Then for all $x\in X$
\[
   \varphi(x) = \int gfx\, d\mu \, ,
\]
hence $\varphi= gf\in X^{\times}$. But $X'= X^{\times}$ implies that $X$ is
$\sigma$-order continuous (see e.g. \cite[II,p.29]{LT79}).
}
\end{rem}

\bn
The dual statement of Corollary~\ref{cor-2} is

\begin{cor}\label{cor-6}
Let $1\le r< \infty$ and $F$ be a Banach space. Then the following are
equivalent:
\begin{itemize}
\item[{\rm (1)}]
Every $r'$-summing operator on $F$ is 1-summing
(=: the identity on $F$ is $(r',1)$-mixing).
\item[{\rm (2)}]
Every operator $T: X\lra F$, $X$ a Banach lattice, is $r$-concave.
\item[{\rm (3)}]
For every operator $T: X(\mu)\lra F$, $X(\mu)$ an $r$-convex quasi
Banach function space, there is some linear $0\le \varphi:
X^r\lra \R$
with $\sup_{\| x\|_{X}\le 1} \varphi \bl |x|^r\br^{1/r} < \infty$
such that for all $x\in X$
\[
   \| Tx \|_F\le \varphi (|x|^r)^{1/r}\, .
\]
\end{itemize}
\end{cor}
Again we only sketch the

\begin{proof}
As mentioned in Remark~\ref{rem-3}, (c), statement (1) by trace duality is
equivalent
to the fact that each operator $T: \e_{\infty}\lra F$ is $r$-summing
(\cite[32.2]{DF93}), and by the definitions an operator $T:
\e_{\infty}\lra F$ is
$r$-summing if and only if it is $r$-concave. On the other hand, by
Krivine's localization technique from \cite{Kr73} (see again the proof of
\cite[II,1.f.14]{LT79}) each $T: \e_{\infty}\lra F$ is $r$-concave if and
only if
(2) holds. This shows that (1) $\car\hspace*{-4.2mm}\cale$
(2). Clearly, (2) implies (3) by Corollary~\ref{cor-5}. In order to
prove the implication (3) $\car$ (2) we only have to check that each $T:
\e_{\infty}\lra F$ is $r$-concave. But since $\e_{\infty}$ is $r$-convex,
this follows by an easy calculation from the inequality in (3).
\end{proof}

\bn
Finally, we mention

\begin{cor}\label{cor-7}
Let $1\le r < \infty$ and $F$ be a Banach lattice. Then the following are
equivalent:
\begin{itemize}
\item[{\rm (1)}]
$F$ is $r$-concave.
\item[{\rm (2)}]
Every positive operator $T: X\lra F$, $X$ a Banach lattice, is
$r$-concave.
\item[{\rm (3)}]
Each positive operator $T: X(\mu)\lra F$, $X(\mu)$ an $r$-convex
Banach function space over $\mu$, satisfies {\rm (3)} of
Corollary~{\rm \ref{cor-6}}.
\end{itemize}
\end{cor}

The implications (1) $\car$ (2) $\car$ (3) follow as in the proof of
Corollary~\ref{cor-2} and by Theorem~\ref{theo-2}; for (3) $\car$ (1)
note first that each positive $T: \e_{\infty}\lra F$ by assumption is
$r$-summing (= $r$-concave) which by Maurey's result already mentioned in
Remark~\ref{rem-3}, (c) gives (1).


\subsection{The Maurey-Rosenthal theorem for operators between
vector-\-valued quasi Banach function spaces}

Finally, we state a variant of Theorem~\ref{theo-2} for vector-valued
quasi
Banach spaces $X(\mu,E)$ as described at the beginning of this section --
the following consequence of Theorem~\ref{theo-2} and the representation
given in (E) includes Corollary~\ref{cor-2} and \ref{cor-5} as
special cases.

\begin{cor}\label{cor-8}
Let $0 < r < \infty$. Assume that $X(\mu,E)$ and $Y(\nu,F)$ are
vector-valued quasi Banach spaces, $X(\mu)$ $r$-convex and $\sigma$-order
continuous, and $Y(\nu)$ $r$-concave. Then for
each homogeneous operator $T: X(\mu,E) \lra Y(\nu,F)$ such that for all
$x_1,\ldots, x_n \in X(\mu,E)$
\begin{equation}\label{eq-4.1}
  \bdv \bl \tsum \bdv Tx_k\bdv^r_F\br^{1/r}\bdv_Y \le
  \bdv \bl \tsum \bdv x_k\bdv^r_E\br^{1/r}\bdv_X
\end{equation}
there are weights
\begin{align*}
    & 0\le \omega_1\in L_0(\mu)\q \text{with}\q \sup_{\| x\|_X\le 1} \,
      \bdv\omega_1^{1/r}x\bdv _{L_r(\mu)} \le M^{(r)}(X)\\[1ex]
    & 0\le \omega_2\in L_0(\nu)\q \text{with}\q \sup_{\| y\|_{L_r(\nu)}\le
        1}\,\bdv\omega_2^{1/r}y\bdv_Y \le M_{(r)}(Y)
\end{align*}
such that for all $x\in X$
\[
   \int \frac{\|Tx\|^r_F}{\omega_2} \, d\nu \le \int \|x\|^r_E \,\omega_1
                    \, d\mu\, .
\]
\end{cor}

Clearly, under appropriate additional assumptions this result can also be
formulated as a factorization theorem.

Let us give an example in the setting of Lorentz spaces: assume that $0 <
p_1 < r < q_1 < \infty$, $0 < p_2 \le r \le q_2 \le \infty$. Then each
homogeneous operator $T: L_{q_1,q_2}(\mu,E) \lra L_{p_1,p_2}(\nu,F)$ ($E$
and $F$ two quasi Banach spaces) which satisfies a vector-valued norm
inequality
\[
  \bdv \bl \tsum \| Tx_k\|^r_F\br^{1/r}\bdv_{L_{p_1,p_2}} \le
  \bdv \bl \tsum \| x_k\|_E^r\br^{1/r}\bdv_{L_{q_1,q_2}}\, ,
\]
satisfies a weighted norm inequality
\[
  \int\frac{\| Tx\|^r_F}{\omega_2}\, d\nu \le \int \| x\|_E^r \, \omega_1\,
   d\mu
\]
(and vice versa).

As already mentioned, in the scalar valued case positive operators $T:
X\lra Y$ between Banach lattices satisfy for all $1\le r < \infty$ a
vector-valued norm inequality
\[
  \bdv \bl \tsum \bv Tx_k\bv^r\br^{1/r}\bdv_Y \le
  \| T \| \bdv \bl \tsum \bv x_k\bv^r\br^{1/r}\bdv_X\, ;
\]
for $r=2$ the famous Krivine-Grothendieck inequality states that this is
even true for all operators (see \cite{Kr73} and \cite[II, 1.f.14]{LT79}).

\begin{rem}\label{rem-5} {\rm
For which pairs $(E,F)$ of Banach spaces is \eqref{eq-4.1}
satisfied for {\em all\/} operators
$T: L_q(\mu,E)\lra L_p(\nu,F)$ acting between spaces of Bochner integrable
functions? For which pairs of Banach lattices $(E,F)$ does this hold for
{\it all\/} positive $T$?

\bsk
In \cite{DJI} it is shown that whenever $E$ or $F$ has the approximation
property (or is $K$-convex or a Banach lattice), then the following are
equivalent:
\begin{itemize}
\item[(1)]
$E'$ and $F$ have cotype 2 .
\item[(2)]
There is $c\ge 0$ such that for all $n$ and all operators $T:
\e^n_\infty(E)\lra \e^n_1(F)$
\[
   \bdv \bl \tsum \| Tx_k\|^2_F\br^{1/2}\bdv_{\e^n_1} \le
     c\, \| T\| \bdv \bl \tsum \bdv
     x_k\bdv^2_E\br^{1/2}\bdv_{\e^n_{\infty}}
   \text{ for all } x_1,\ldots, x_m\in \e^n_\infty(E) \, .
\]
\end{itemize}
}
\end{rem}

Moreover, we know from a monotonicity argument in \cite{DJI} that (1)
implies (2) if one replaces
$\e^n_1(F)$ by $\e^n_p$(F) and $\e^n_{\infty}(E)$ by $\e^n_q(E)$
$(1\le p,q\le \infty)$; using
conditional expectation operators, it hence can be checked that if $E$ or
$F$ has the approximation property (or is $K$-convex or a Banach lattice)
and $E',F$ have
cotype 2, then there is a constant $c\ge 0$ such that for all $1\le p,q<
\infty$ and all operators $T: L_q(\mu,E)\lra L_p(\nu,F)$
\[
  \bdv \bl \tsum \| Tx_k\|^2_F\br^{1/2}\bdv_{L_p(\nu)} \le
  c\, \|T\| \bdv \bl \tsum \| x_k\|^2_E\br^{1/2}\bdv_{L_q(\mu)}
  \q\text{for all}\q x_1,\ldots, x_n\in L_q(\mu,E)\,.
\]
 Hence, if moreover
$1\le p< r = 2 < q < \infty$, then all such $T$ allow a factorization as in
Corollary~\ref{cor-5}.

\bsk
A result which will be shown in \cite{DJII} states that for two Banach
lattices $E,F$ and $1\le p< r< q<\infty$ all positive operators $T:
L_q(\mu,E)\lra L_p(\nu,F)$ satisfy \eqref{eq-4.1} (up to some constant) if
and only if $E$ is $r$-convex and $F$ $r$-concave.


\end{document}